\documentclass[11pt,notitlepage,twoside]{article}

\usepackage{latexsym}
\usepackage{amsmath}
\usepackage{amsfonts}
\usepackage{amssymb}

\newcommand{\be}{\begin{equation}}
\newcommand{\ee}{\end{equation}}
\newcommand{\hf}{\hfill $\diamondsuit$}
\newenvironment{pf}{\noindent{\bf Proof.}\enspace}{
\hfill Q.E.D.}

\usepackage{latexsym}
\newtheorem{thm}{Theorem}[section]
\newtheorem{pro}{Proposition}[section]
\newtheorem{lem}{Lemma}[section]
\newtheorem{rem}{Remark}[section]

\newtheorem{df}{Definition}[section]
\newtheorem{ex}{Examples}[section]
\numberwithin{equation}{section}
\DeclareFontFamily{T1}{cmr}{\hyphenchar\font=-1}
\begin{document}
\centerline {\bf \Large On a Riesz Basis of   }
\vskip 0.1 cm
\centerline {\bf \Large Diagonally Generalized Subordinate   Operator Matrices  }
\vskip 0.1 cm
\centerline {\bf \Large  and Application to a Gribov  Operator Matrix}
\centerline {\bf \Large  in Bargmann Space}
\vskip 1 cm

{\centerline{\bf{B. Abdelmoumen, A. Damergi  and Y. Krichene}}}
\vskip 0.1cm
Department of Mathematics, Faculty of Sciences of Sfax, University of Sfax, Soukra\vskip 0.2 cm \centerline{Road Km 3.5,
B. P. 1171, 3000 Sfax, Tunisia.}

\vskip 0.2 cm

{\centerline{\small {e-mail:  boulbeba.abbelmoumen@ipeis.rnu.tn}}}
{\centerline{\small {e-mail: damergialaeddine1981@gmail.com}}}
{\centerline{\small {e-mail: yousraakrichene@gmail.com}}}

\vspace{3mm}
\normalsize
\vspace{3mm} \hspace{.05in}\parbox{4.5in}{{\bf\small Abstract. }{\small  \sloppy{

In this paper, we study  the change of  spectrum and   the existence  of Riesz  bases  of   specific classes  of  $n\times n$ unbounded operator matrices,  called:  diagonally  and off-diagonally  generalized subordinate block operator matrices. An application    to a  $n\times n$ Gribov   operator matrix acting on a sum of  Bargmann   spaces, illustrates  the abstract  results. As example,  we consider a particular Gribov operator matrix by taking special   values of the real  parameters   of Pomeron.}}}
\section{ Introduction}
Block operator matrices are matrices with entries that are linear operators between Banach or Hilbert spaces.  They play a crucial role    in several  active fields of research: in control and systems theory, in differential equations, in evolution problems and   in   mathematical physics with its various   applications, such as  hydrodynamics, magnetohydrodynamics,  elasticity theory, optimal control,  mechanics and quantum  mechanics.  The study of spectral proprieties  of this  type of operators is then of great importance, but also of complexity,   and the situation becomes more complicated  for $n\geq3$. In effect, most of   the works that have been carried out  on this subject have been done  in the particular  case $n=2$, and  we can quote  for instance
\cite{A.B9, ACM, ALMS, Engel, Engel1, LT, HLT, Nagel1, Nagel2,  Tretter1}.  Among the  works that have dealt with  the general case $n\geq2$, we can cite  \cite{Tretter3} and \cite{TW}.  In this latter, Tretter et al. introduced the block numerical range of bounded  $n\times n$ operator matrices  and proved the spectral inclusion property. In \cite{Tretter3},  the results obtained in   \cite{ Tretter1} and \cite{TW}  was generalized to unbounded \emph{ diagonally}  (respectively \emph{off-diagonally}) \emph{ dominant $n\times n$ operator matrices}, which  are matrices such that the off-diagonal entries of each column are relatively bounded to  the diagonal entry of the same column (respectively  the diagonal entry of each column is relatively bounded to  each off-diagonal entry of the same column).
\vskip0.3 cm

 In this paper, we introduce   specific   matrices of the last classes  called \emph{diagonally}  and  \emph{off-diagonally}  \emph{generalized subordinate matrices}, for which  we will  characterize the location of the spectrum and  study the existence of Riesz bases of generalized eigenvectors.

 The results of this paper, are of importance for application to a $n\times n$  Gribov     operator matrix, which is an operator matrix whose entries are Gribov operators.

Notice that the Gribov operators form  a family of non-self-adjoint operators that govern the   Reggeon field theory. This theory was introduced   by  Gribov  \cite{Gribov} in 1967, in order  to study strong interactions, i.e., the interaction between protons and neutrons among other less stable particles. These operators  was originally    constructed as polynomials in   the standard creation and   annihilation operators. A representation   space of this family is the Bargmann space of    analytic  functions on $\mathbb{C}^{d}$. In zero transverse dimension ($d=1$),  the Bargmann space is the following:
$$\hfill \mathcal{B}=\left\{ \varphi:\mathbb{C}\longrightarrow \mathbb{C}
\hbox{~entire such that~}\displaystyle \int_{\mathbb{C}}e^{-\vert
z\vert^2}\vert \varphi(z)\vert^2
dx dy<\infty\hbox{~and~}\varphi(0)=0\right\}.$$
\noindent It forms  with  the inner product
 $$\langle\varphi, \psi\rangle \longrightarrow \frac{1}{\pi}\int_{\mathbb{C}}e^{-|z|^{2}}\varphi(z)\overline{\psi}(z)dx dy,$$
\noindent   a Hilbert space. An orthonormal basis for $\mathcal{B}$, is $\displaystyle\{z\mapsto e_{k}(z)=\frac{z^{k}}{\sqrt{k!}}\}_{k\in \mathbb{N}}$.
   \vskip0.2 cm

 In this representation, the standard  annihilation operator $A$  and the standard creation  operator $A^{*}$, are defined by:
\par \vskip0.3cm$\hfill\left\{
\begin{array}{l}
A:{\mathcal{D}}(A)\subset \mathcal{B}\longrightarrow
\mathcal{B}\\
\displaystyle{\quad\quad\quad\quad\quad\,\,\,\,\,\,\!\!\!\varphi
\longrightarrow A\varphi=\frac{d\varphi}{dz}}\\\\
{\mathcal{D}}(A)=\{\varphi\in \mathcal{B} \,\,\hbox{such that}\,\,A\varphi\in \mathcal{B}\}
\end{array}
\right.\hfill$\par \vskip0.3cm\noindent and\par \vskip0.3cm

$\hfill\left\{
\begin{array}{l}
A^{\ast}:{\mathcal{D}}(A^{\ast})\subset \mathcal{B}\longrightarrow
\mathcal{B}\\
\displaystyle{\quad\quad\quad\quad\quad\quad\,\,\,\,\!\!\!\varphi
\longrightarrow A^{\ast}\varphi:z\mapsto z\varphi(z)}\\\\
{\mathcal{D}}(A^{\ast})=\{\varphi\in \mathcal{B}\,\,\hbox{such that}\,\,A^{\ast}\varphi\in \mathcal{B}\}.
\end{array}
\right.\hfill$\par \vskip0.3cm

The following operators play a crucial role in Reggeon field theory.
$$\left\{
\begin{array}{l}
G:{\mathcal{D}}(G)\subset \mathcal{B}\longrightarrow
\mathcal{B}\\
\displaystyle{\quad\quad\quad\quad\quad\quad\,\,\,\,\!\!\!\varphi
\longrightarrow G\varphi=A^{{*}3}A^3\varphi}\\\\
{\mathcal{D}}(G)=\{\varphi\in \mathcal{B}\,\,\hbox{such that}\,\,G\varphi\in \mathcal{B}\},
\end{array}
\right.$$

$$\left\{
\begin{array}{l}
S:{\mathcal{D}}(S)\subset \mathcal{B}\longrightarrow
\mathcal{B}\\
\displaystyle{\quad\quad\quad\quad\quad\quad\,\,\,\,\!\!\!\varphi
\longrightarrow S\varphi=A^{{*}2}A^2\varphi}\\\\
{\mathcal{D}}(S)=\{\varphi\in \mathcal{B}\,\,\hbox{such that}\,\,S\varphi\in \mathcal{B}\},
\end{array}
\right.$$
$$\left\{
\begin{array}{l}
H_{0}:{\mathcal{D}}(H_{0})\subset \mathcal{B}\longrightarrow
\mathcal{B}\\
\displaystyle{\quad\quad\quad\quad\quad\quad\,\,\,\,\!\!\!\varphi
\longrightarrow H_{0}\varphi=A^{{*}}A\varphi}\\\\
{\mathcal{D}}(H_{0})=\{\varphi\in \mathcal{B}\,\,\hbox{such that}\,\,H_{0}\varphi\in \mathcal{B}\}
\end{array}
\right.$$
and
$$\left\{
\begin{array}{l}
H_{1}:{\mathcal{D}}(H_{1})\subset \mathcal{B}\longrightarrow
\mathcal{B}\\
\displaystyle{\quad\quad\quad\quad\quad\quad\,\,\,\,\!\!\!\varphi
\longrightarrow H_{1}\varphi=A^{\ast}(A^{\ast}+A)A\varphi}\\\\
{\mathcal{D}}(H_{1})=\{\varphi\in \mathcal{B}\,\,\hbox{such that}\,\,H_{1}\varphi\in \mathcal{B}\}.
\end{array}
\right.$$

 \vskip0.2 cm  En effect, a  third degree representative of this theory, is the following Gribov operator:
 \begin{eqnarray}\label{GBO3}
H_{\lambda'',\,\lambda',\,\mu,\,\lambda} & = &\lambda''G+\lambda'S+\mu H_{0} +i \lambda H_{1},
\end{eqnarray}
where $\textbf{i}^{2} =-1$, the real parameters   $\lambda''$, $\lambda'$ and  $\lambda$   are respectively  the magic, the four and   the triple coupling of
Pomeron and the real  $\mu$ is the Pomeron intercept.

\vskip 0.3 cm    The  paper is organized as follows:  In Section
$2$, we study the location   of spectrum and the  existence of Riesz bases of   diagonally and off-diagonally generalized subordinate $n\times n$  operator matrices.
In section  $3$,  the main results are illustrated by an application  to a  problem of a $n\times n$  Gribov  operator matrix acting on  a sum of  $n$ identical  Bargmann spaces.
 \section{ Main results}
Let $\mathcal{H}_{1},\ldots,\mathcal{H}_{n}$ be Hilbert spaces.  In the Hilbert space   $\mathcal{H}:=\displaystyle\bigoplus_{j=1}^{n}\mathcal{H}_{j}$,  we consider    a $n\times n $  operator matrix  $\mathcal{M}=(A_{ij})_{1\leq i,j\leq n}$,
 where the entries $A_{ij}$ are unbounded linear operators from $\mathcal{D}(A_{ij}) \subset \mathcal{H}_{j}$ to $\mathcal{H}_{i}$. Obviously, $\mathcal{M}$ has as natural domain
$$\mathcal{D}(\mathcal{M}):=\displaystyle\bigoplus_{j=1}^{n}\bigcap_{i=1}^{n}\mathcal{D}(A_{ij})$$ and  if the entries are densely defined, then $\mathcal{M}$ is too.
\vskip0.2 cm

 For the sequel, we adopt  the following  decomposition:  $$\mathcal{M}=\mathfrak{D}+\mathcal{R}$$ with $$\mathfrak{D}=diag(A_{11},...,A_{nn})\,\,\hbox{and}\,\,
\mathcal{R}=((1-\delta_{ij})A_{ij})_{1\leq i, j\leq n},$$
where $\delta_{ij}$ denote the Kronecker delta of $(i,j)$.  For the sake of simplicity, we denote by  $\Omega_{n}$, the following set:
$$\Omega_{n}:=\{(i,j)\in\{1,\ldots,n\}^{2}\,\,\hbox{such\,\,that}\,\,i\neq j\}.$$

\vskip0.4 cm
 The new concepts of diagonally  and off-diagonally generalized subordinate  operator matrices are given in the following definition.
\begin{df} The  operator matrix $\mathcal{M}$ is called
\begin{enumerate}
  \item   diagonally generalized subordinate,  if for all $(i,j)\in \Omega_{n}$,

   \hfill  $A_{ij}$   is  generalized subordinate to  $A_{jj}, \hfill$

     i.e, there exist $N^{(i,j)}\in \mathbb{N}^{*},$  $\displaystyle\{p_{k}^{(i,j)}\}_{k=1}^{N^{(i,j)}}\subset [0,1]$   and   $\displaystyle\{b_{k}^{(i,j)}\}_{k=1}^{N^{(i,j)}}\subset \mathbb{R}_{+}$, such that

$\hfill\|A_{ij}U_j\|\leq  \displaystyle\sum_{k=1}^{N^{(i,j)}} \displaystyle b_{k}^{(i,j)}\|A_{jj}U_j\|^{p_{k}^{(i,j)}}\|U_j\|^{1-p_{k}^{(i,j)}},\,\,\forall \,\,U_j\in \mathcal{D}(A_{jj}).\hfill$
  \item off-diagonally generalized subordinate, if for all $(i,j)\in \Omega_{n}$,

 \hfill   $A_{jj}$  is generalized subordinate to  $A_{ij},\hfill$

   i.e, there exist $N^{'(i,j)}\in \mathbb{N}^{*},$  $\{p_{k}^{'(i,j)}\}_{k=1}^{N^{'(i,j)}}\subset [0,1]$ and $\displaystyle\{b_{k}^{'(i,j)}\}_{k=1}^{N^{'(i,j)}}\subset \mathbb{R}_{+}$, such that

$\hfill\|A_{jj}U_j\|\leq   \displaystyle\sum_{k=1}^{N^{'(i,j)}}  \displaystyle b_{k}^{'(i,j)}\|A_{ij}U_j\|^{p_{k}^{'(i,j)}}\|U_j\|^{1-p_{k}^{'(i,j)}}, \forall \,\,U_j\in \mathcal{D}(A_{ij}).\hfill$\hf
\end{enumerate}
\end{df}
\vskip0.3 cm
\begin{rem}
Note that the concept of \emph{generalized subordinate perturbations} was  introduced  in \cite{AL} as a  natural  generalization  of   the  notion  of  \emph{$p$-subordinate perturbations}  studied by Krein (see \cite{S.G}),  Markus (see\cite{M1})  and Wyss   (see \cite{w1}  and  \cite{w}).\hf
\end{rem}
\vskip0.3 cm
Now, we will state    a crucial result   that we will use in the rest of this paper.
\begin{lem} \label{l0}
The following hold:
\begin{description}
  \item[(i)] Suppose   $\mathcal{M}$ is  diagonally generalized subordinate, then $\mathcal{R}$   is generalized subordinate to $\mathfrak{D}$, i.e,
  $\mathcal{D}(\mathfrak{D})\subset\mathcal{D}(\mathcal{R})$ and
   $\exists N\in \mathbb{N}^{*},\,\, \displaystyle\{p_{k}\}_{k=1}^{N}\subset [0,1]$   and  $\displaystyle \{b_{k}\}_{k=1}^{ N}\subset \mathbb{R}_{+}$, such that

    $\hfill\|\mathcal{R}U\|\leq \displaystyle\sum_{k=1}^{N}b_{k}\|\mathfrak{D}U\|^{p_{k}}\|U\|^{1-p_{k}},\ \forall U\in \mathcal{D}(\mathfrak{D}).\hfill$
  \item[(ii)] Suppose that $\mathcal{M}$ is  off-diagonally generalized subordinate, then  $\mathfrak{D}$  is generalized subordinate to $\mathcal{R}$, i.e,
   $\mathcal{D}(\mathcal{R})\subset\mathcal{D}(\mathfrak{D})$ and  $\exists N'\in \mathbb{N}^{*},\,\,\{p^{'}_{k}\}_{k=1}^{ N'}\subset [0,1]$ and $\{b^{'}_{k}\}_{k=1}^{ N'}\subset \mathbb{R}_{+}$, such that

$\hfill\|\mathfrak{D}U\|\leq \displaystyle\sum_{k=1}^{N'}b^{'}_{k}\|\mathcal{R}U\|^{p^{'}_{k}}\|U\|^{1-p^{'}_{k}},\ \forall U\in \mathcal{D}(\mathcal{R}).\hfill$\hf
\end{description}
\end{lem}
\begin{pf}
We prove \textbf{(i)}, the proof of \textbf{(ii)} is analogue. Let $ U= (U_1,...,U_n)\in \displaystyle \bigoplus_{j=1}^{n}\mathcal{D}(A_{jj})$,
\begin{eqnarray*}
 \|\mathcal{R}U\| &=& \sum_{i=1}^n\|\sum_{j=1\atop j\neq i}^{n}A_{ij}U_j\|\\
 &\leq&  \sum_{(i,j)\in \Omega_{n}}\|A_{ij}U_j\|\\
&\leq& \sum_{(i,j)\in \Omega_{n}}\sum_{k=1}^{N^{(i,j)}}b_{k}^{(i,j)}\|A_{jj}U_j\|^{p_{k}^{(i,j)}}\|U_j\|^{1-p_{k}^{(i,j)}}\\
&\leq& \sum_{(i,j)\in \Omega_{n}}\sum_{k=1}^{N^{(i,j)}}b_{k}^{(i,j)}\|\mathfrak{D} U\|^{p_{k}^{(i,j)}}\|U\|^{1-p_{k}^{(i,j)}}.
\end{eqnarray*}
\noindent   It suffices to take\vskip 0.2 cm
\begin{equation}\label{sys}
    \left\{
  \begin{array}{ll}
    N=\displaystyle\sum_{(i,j)\in \Omega_{n}}N^{(i,j)}, & \hbox{} \\
    \{b_{k}\}_{k=1}^{N}=\displaystyle\bigcup_{(i,j)\in \Omega_{n}}\{b_{k}^{(i,j)} \}_{k=1}^{N^{(i,j)}}, & \hbox{} \\
    \{p_{k}\}_{k=1}^{N}=\displaystyle\bigcup_{(i,j)\in \Omega_{n}}\{p_{k}^{(i,j)} \}_{k=1}^{N^{(i,j)}}. & \hbox{}
  \end{array}\right.
\end{equation}
\noindent  This ends the proof.
\end{pf}
\vskip 0.4 cm

The following theorem  provides a detailed  description of the changed   spectrum of diagonally  generalized  subordinate  $n\times n$ block operator matrices.  For off-diagonally generalized  subordinate  $n\times n$ block operator matrices, the result is similar.
\begin{thm}\label{t4} Suppose  that the  operator matrix  $\mathcal{M}$ is diagonally  generalized  subordinate, with  $\displaystyle\{p_{k}\}_{1\leq k\leq N}\subset[0,1[$. Suppose moreover, that  for all    $j=1,\ldots,n$, $ A_{jj}$  is normal with compact resolvent and  its spectrum  lies on finitely many rays:
$$ \sigma(A_{jj})\displaystyle\subset \bigcup_{l=1}^{m_{j}}e^{i\theta_{jl}}\mathbb{R}_{\geq0};\,\,0\leq\theta_{jl}\leq 2\pi,\,\,l=1,\cdots,m_{j}.$$ Then $\mathcal{M}$ is with compact resolvent and there exists $m\in \mathbb{N}^{*}$, such that
for every $\alpha>\displaystyle\sum_{k=1}^{N}b_{k}$,
there exists $r_0> 0$ satisfying

\hfill$\displaystyle{\sigma(\mathcal{M})\subset B_{r_0}(0)\cup\bigcup_{l=1}^{m}\{  \displaystyle e^{\displaystyle\textbf{i}\theta_l}(x+\textbf{i}y);\  x\geq0\mbox{ and } |y|\leq \max_{1\leq k\leq N}\alpha x^{p_{k}}\}.}$\hfill\hf
\end{thm}
\begin{pf}
 The entries   $ \displaystyle (A_{jj})_{1\leq j\leq n}$  are  normals with compact resolvents. So,  $\mathfrak{D}$ is  too and
$$\sigma(\mathfrak{D})=\bigcup_{s=1}^n\sigma(A_{ss})\subset \displaystyle\bigcup_{s=1}^{n}\bigcup_{l=1}^{m_{s}}\displaystyle e^{\displaystyle\textbf{i}\theta_{sl}}\mathbb{R}_{\geq0}.$$
\noindent Apply  Lemma \ref{l0} and  \cite[ Theorem 2.1]{AL} to the decomposition $\mathcal{M}=\mathfrak{D}+\mathcal{R}$, we obtain the required result.
\end{pf}\vskip 0.3 cm

Since for non-normal operators there is no analogue of the spectral theorem, the existence of Riesz bases is of great importance. Recently, some progress has been made for the existence of the  Riesz bases existence of eigenvectors of $2\times2$ block operators matrices,
(see \cite{ CDJ, Jacob, Kuiper, w1, w}) and $3\times3$ block operators matrices (see \cite{AL}).  For $n>3$, up until now, there have been no results on the existence of Riesz bases for   block operator matrices.  The following theorem  provides different necessary conditions in terms of the spectrum in order to prove the existence of Riesz bases with parentheses for diagonally  generalized  subordinate $n\times n$ block operator matrices.
\begin{thm}\label{t5}
Under the  hypotheses  of Theorem \ref{t4}, suppose moreover that
$$\forall\,s\in\{1,\ldots,n\},\,\,\displaystyle\liminf_{r\rightarrow+\infty}\frac{N(r,A_{ss})}{r^{1-p}}<+\infty,\ \mbox{where }p= \max_{\ 1\leq k\leq N}\{p_{k}\}.$$
Then  the  operator matrix $\mathcal{M}$ admits a Riesz basis with parentheses.\hf
\end{thm}
\begin{pf}
 We have $ \displaystyle N(r,\mathfrak{D})=\sum_{s=1}^n \displaystyle N(r,A_{ss}).$  Hence, the result follows from    Lemma \ref{l0} and \cite[ Theorem 2.2]{AL}.
\end{pf}
\section{ Gribov  Operator  Matrix in Bargmann space}
The Gribov operator matrices are  matrices  whose entries are Gribov operators  between subspaces of the  Bargmann space $\mathcal{B}$. To illustrate the applicability  of abstract results described above, we give an application to a representative of the $n\times n$  Gribov operator matrices in zero transverse dimension. To define this representative, we must  start by giving some  fundamental  results  and introducing  some basic  notations.
 \vskip0.3 cm
We first begin by recalling basic properties of the operators $G$, $S$  and $ H_{0}$, that follow from \cite[Lemma 3, p 112]{Intissar} and
 \cite[Proposition 4, p 112]{Intissar}.
\begin{pro}\label{s13}\par \vskip0.3cm\begin{description}
  \item[(i)] The operators $G$, $S$  and $ H_{0}$ are self-adjoint and with compact resolvents.
  \item[(ii)] $\{e_{k}:\,=\frac{z^{k}}{\sqrt{k!}}\}_{k\geq1}$ is a system of eigenvectors for $ H_{0}$ associated with the eigenvalues $\{n\}_{n\geq1}$ and
 its spectral decomposition  is the following:
 $$H_0\varphi=\sum_{k=1}^{\infty}k\langle\varphi,e_{k}\rangle e_{k},\,\,\forall\,\,\varphi\in {\mathcal{D}}(H_0).$$
 \item[(iii)] $\{e_{k}:\,=\frac{z^{k}}{\sqrt{k!}}\}_{k\geq1}$ is a system of eigenvectors for $ S$ associated with the eigenvalues $\{k(k-1)\}_{k\geq1}$
and  its spectral decomposition  is the following:
 $$S\varphi=\sum_{k=2}^{\infty}k(k-1)\langle\varphi,e_{k}\rangle e_{k},\,\,\forall\,\,\varphi\in {\mathcal{D}}(S).$$
\item[(iv)]  $\{e_{k}:\,=\frac{z^{k}}{\sqrt{k!}}\}_{k\geq2}$ is a system of eigenvectors for $G$ associated with the eigenvalues\par \vskip0.2cm
$\{k(k-1)(k-2)\}_{k\geq2}$ and
 its spectral decomposition  is the following:

 $\hfill G\varphi=\displaystyle\sum_{k=2}^{\infty}(k-2)(k-1)k\langle\varphi,e_{k}\rangle e_{k},\,\,\forall\,\,\varphi\in {\mathcal{D}}(G).\hfill$ $\diamondsuit$\end{description}
\end{pro}\par \vskip0.3 cm

Now, for   $\lambda''\in \mathbb{R}^{*} $,   $(\lambda',\lambda,\mu)\in\mathbb{R}^{3}$ and  $\beta>0$, we introduce the following operators:
$$\left\{
\begin{array}{l}
G_{\lambda''}:{\mathcal{D}}(G_{\lambda''})\subset \mathcal{B}\longrightarrow
\mathcal{B}\\\\
\quad\quad\quad\quad\quad\quad\,\,\,\,\!\!\!\varphi
\longrightarrow G_{\lambda''}\varphi=\lambda''G\varphi\\\\
{\mathcal{D}}(G_{\lambda''})={\mathcal{D}}(G),
\end{array}
\right.$$
\vskip0.3 cm

$$\left\{
\begin{array}{l}
H_{0}^{\beta}:{\mathcal{D}}(H_{0}^{\beta})\subset \mathcal{B}\longrightarrow
\mathcal{B}\\
\displaystyle{\quad\quad\quad\quad\quad\quad\,\,\,\,\!\!\!\varphi
\longrightarrow H_{0}^{\beta}\varphi:=\sum_{k=1}^{\infty}k^{\beta}\langle\varphi,e_{k}\rangle e_{k}}\\\\
{\mathcal{D}}(H_{0}^{\beta})=\left\{\varphi\in \mathcal{B} \,\,\hbox{such that}\,\,{\displaystyle\sum_{k=1}^{\infty}k^{2\beta}|\langle\varphi,e_{k}\rangle|^{2}<\infty}\right\},
\end{array}
\right.$$ \vskip0.3 cm
\noindent and \vskip0.3 cm
$$\left\{
\begin{array}{l}
H^{\beta}_{\lambda',\lambda,\mu}:{\mathcal{D}}(H^{\beta}_{\lambda',\lambda,\mu})\subset \mathcal{B}\longrightarrow
\mathcal{B}\\
\quad\quad\quad\quad\quad\quad\,\,\,\,\!\!\!\varphi
\longrightarrow H^{\beta}_{\lambda',\lambda,\mu}\varphi=[\lambda'S+i \lambda H_{1} +\mu  H_{0}^{\beta}]\varphi\\\\
{\mathcal{D}}(H^{\beta}_{\lambda',\lambda,\mu})=\left\{\varphi\in \mathcal{B} \,\,\hbox{such that}\,\,H^{\beta}_{\lambda',\lambda,\mu}\varphi\in \mathcal{B}\right\}.
\end{array}
\right.$$
 \vskip0.4 cm

The  following  lemma   provides   several useful inequalities for the  operators we defined previously.
\begin{lem}\label{lemH1}
 Put  $c_{1}=5$, $c_{2}=\displaystyle\sqrt{1+2^{6}}$ and $c_{3}=1+2\sqrt{2}$.  The following inequalities hold trues:
\begin{description}
\item[(i)] $\displaystyle\|H_{0}^{3}\varphi\|\leq c_{1} \|G\varphi\|+c_{2}\|\varphi\|$, $\forall\varphi\in\mathcal{D}(G)$.\vskip 0.3cm
    \item[(ii)]  $\displaystyle\|H_{0}^{\beta}\varphi\|\leq \displaystyle\|H_{0}^{3}\varphi\|^{\frac{\beta}{3}}\|\varphi\|^{1-\frac{\beta}{3}}\leq \displaystyle c_{1}^{\frac{\beta}{3}} \|G\varphi\|^{\frac{\beta}{3}}\|\varphi\|^{1-\frac{\beta}{3}}+c_{2}^{\frac{\beta}{3}}\|\varphi\|,$  $\forall\beta\in]0,3[$ and
    $\forall\varphi\in\mathcal{D}(H_{0}^{\beta})$.\vskip 0.3cm
  \item[(iii)] $\displaystyle\|S\varphi\|\leq \displaystyle\|H_{0}^{3}\varphi\|^{\frac{2}{3}}\|\varphi\|^{\frac{1}{3}}\leq \displaystyle c_{1}^{\frac{2}{3}} \|G\varphi\|^{\frac{2}{3}}\|\varphi\|^{\frac{1}{3}}+c_{2}^{\frac{2}{3}}\|\varphi\|,$
    $\forall\varphi\in\mathcal{D}(S)$.\vskip 0.3cm
\item[(iv)] $\displaystyle\|H_{1}\varphi\|\leq c_{3}\displaystyle\|H_{0}^{\frac{3}{2}}\varphi\|\leq \displaystyle c_{3} \displaystyle\|H_{0}^{3}\varphi\|^{\frac{1}{2}}\|\varphi\|^{\frac{1}{2}}\leq \displaystyle c_{3} c_{1}^{\frac{1}{2}} \|G\varphi\|^{\frac{1}{2}}\|\varphi\|^{\frac{1}{2}}+c_{3} c_{2}^{\frac{1}{2}}\|\varphi\|, $
    $\forall\varphi\in\mathcal{D}(H_{1})$.\hf
\end{description}
 \end{lem}\vskip 0.3 cm
 \begin{pf}
 Firstly, for all $(a,b)\in \mathbb{R}_{+}^{2}\backslash\{(0,0)\}$, we have the following  inequalities:
\begin{eqnarray}
    (a+b)^{\alpha}\leq a^{\alpha}+b^{\alpha},&& \hbox{if}\,\,0\leq\alpha\leq 1,\label{Eq1} \\
   a^{\alpha}+b^{\alpha}\leq (a+b)^{\alpha}, &&\hbox{if}\,\,\alpha\geq 1. \label{Eq2}
\end{eqnarray}

\noindent \textbf{(i)} Let $\varphi\in\mathcal{D}(G).$ By using the spectral decomposition of $H_{0}^{3}$ and $G$ and the fact that $k^{2}\leq 5(k-2)(k-1)$ for all $k\geq3$, we have
 \begin{eqnarray*}
   \sum_{k=1}^{\infty}k^{6}\mid \langle\varphi,e_{k}\rangle\mid^{2} &\leq&\mid \langle \varphi,e_{1}\rangle\mid^{2}+2^{6}\mid \langle \varphi,e_{2}\rangle\mid^{2} +\sum_{k=3}^{\infty}[5(k-2)(k-1)k]^{2}\mid \langle \varphi,e_{k}\rangle\mid^{2} \\
    &\leq&  (1+2^{6})\|\varphi\|^{2}+25\|G\varphi\|^{2}.
 \end{eqnarray*}
 Apply Inequality \eqref{Eq1}, the inequality  follows.
\vskip 0.3cm
 \noindent \textbf{(ii)} Let $\varphi\in \mathcal{D}(H_{0}^{\beta})$. The first inequality follows  from the application of the H\^{o}lder inequality:
$$ \sum_{k=1}^{\infty}k^{2\beta}|\langle\varphi,e_{k}\rangle|^{2} \leq
\Big(\sum_{k=1}^{\infty}k^{6}|\langle\varphi,e_{k}\rangle|^{2}\Big)^{\frac{\beta}{3}}
\Big(\sum_{k=1}^{\infty}k^{2\beta}|\langle\varphi,e_{k}\rangle|^{2}\Big)^{1-\frac{\beta}{3}}.$$
For the second inequality, it suffices to apply \textbf{(i)} and inequality \eqref{Eq1}.\vskip 0.3cm
\noindent Concerning the first two inequalities of \textbf{(iii)} and \noindent \textbf{(iv)}, one can see \cite[Lemme 4.1]{CDJ}. The other two inequalities, result from $\textbf{(ii)}$, \eqref{Eq1} and \eqref{Eq2}.
\end{pf}
\vskip 0.4 cm
The next proposition follows from   Lemma \ref{lemH1}.
\begin{pro}\label{ff} Let $(\lambda'',(\lambda',\lambda,\mu),\beta)\in \mathbb{R}^{*}\times\mathbb{R}^{3}\times]0,3[$.\\
 $H^{\beta}_{\lambda',\mu,\lambda}$ is $\{\frac{\beta}{3},\frac{1}{2},\frac{2}{3}\}$-generalized subordinate to $G_{\lambda''}$ with bound
 $\Big\{ \frac{\mid\mu\mid c_{1}^{\frac{\beta}{3}}}{\mid\lambda''\mid^{\frac{\beta}{3}}}, \frac{\mid\lambda\mid c_{3} c_{1}^{\frac{1}{2}}}{\mid\lambda''\mid^{\frac{1}{2}}},   \frac{\mid\lambda'\mid c_{1}^{\frac{2}{3}}}{\mid\lambda''\mid^{\frac{2}{3}}}\Big\}$.  \hf\end{pro}
\begin{pf}
 \begin{eqnarray*}
 \| H^{\beta}_{\lambda',\lambda,\mu}\varphi\| & \leq & \mid\mu\mid \|H^{\beta}_{0}\varphi\| + \mid\lambda\mid \|H_{1}\varphi\|+ \mid\lambda'\mid \|S\varphi\| \\
   &\leq&
   \displaystyle
   \frac{\mid\mu\mid c_{1}^{\frac{\beta}{3}}}{\mid\lambda''\mid^{\frac{\beta}{3}}}\|G_{\lambda''}\varphi\|^{\frac{\beta}{3}}\|\varphi\|^{1-\frac{\beta}{3}}+ \displaystyle \frac{\mid\lambda\mid c_{3} c_{1}^{\frac{1}{2}}}{\mid\lambda''\mid^{\frac{1}{2}}}  \|G_{\lambda''}\varphi\|^{\frac{1}{2}}\|\varphi\|^{\frac{1}{2}}+ \displaystyle  \frac{\mid\lambda'\mid c_{1}^{\frac{2}{3}}}{\mid\lambda''\mid^{\frac{2}{3}}}\|G_{\lambda''}\varphi\|^{\frac{2}{3}}\|\varphi\|^{\frac{1}{3}}
 \end{eqnarray*}
 $\hfill            +c_{2}^{\frac{1}{3}}\mid\mu\mid\|\varphi\|+\displaystyle  c_{3} c_{2}^{\frac{1}{2}}\mid\lambda\mid\|\varphi\|+ \mid\lambda'\mid c_{2}^{\frac{2}{3}}\|\varphi\|.                 \hfill$\end{pf}
\vskip 0.5 cm

Now, for $(\lambda''_{j})_{1\leq j\leq n}\subset\mathbb{R}^{*} $ and $\big((\lambda'_{ij},\lambda_{ij},\mu_{ij}),\beta_{ij})\big)_{(i,j)\in \Omega_{n}} \subset \mathbb{R}^{3}\times]0,3[$,  we   consider the $n\times n$ Gribov  operator  matrix $\mathcal{M}^{\beta}_{\lambda'',\lambda',\lambda,\mu}$ defined by:
\begin{equation}\label{GMAT}
   \mathcal{M}^{\beta}_{\lambda'',\lambda',\lambda,\mu}:=\mathfrak{D}_{\lambda''}+\mathcal{R}^{\beta}_{\lambda',\lambda,\mu},
\end{equation}
where
$$\mathfrak{D}_{\lambda''}=Dig(G_{\lambda''_{1}},\ldots,G_{\lambda''_{n}})\,\,\hbox{and}\,\, R^{\beta}_{\lambda',\lambda,\mu}=\big( (1-\delta_{ij})H^{\beta_{ij}}_{\lambda'_{ij},\lambda_{ij},\mu_{ij}}\big)_{1\leq i,j\leq n}.$$
\vskip 0.4 cm

\begin{rem} The Gribov  operator  matrix $\mathcal{M}^{\beta}_{\lambda'',\lambda',\lambda,\mu}$ can be written as follows:
  \begin{equation}\label{Gdic}
    \mathcal{M}^{\beta}_{\lambda'',\lambda',\lambda,\mu}:=\mathfrak{D}_{\lambda''}+\mathfrak{S}_{\lambda'}+\mathfrak{H}_{\lambda}+\mathfrak{H}_{\mu}^{\beta},
  \end{equation}
where
 $$\mathfrak{S}_{\lambda'}=\mathcal{R}^{\beta}_{\lambda',0,0},\,\,\mathfrak{H}_{\lambda}=\mathcal{R}^{\beta}_{0,\lambda,0}\,\,\hbox{and}\,\,\mathfrak{H}_{\mu}^{\beta}= \mathcal{R}^{\beta}_{0,0,\mu}.$$
 According to this decomposition, $\mathcal{M}^{\beta}_{\lambda'',\lambda',\lambda,\mu}$ can be considered as a representative for the $n\times n$ Gribov operator matrices. \hf
\end{rem}
\vskip 0.3 cm

  For the sequel, we put  $$ b_{1}^{(i,j)}=\frac{\mid\mu_{ij}\mid c_{1}^{\frac{\beta_{ij}}{3}}}{\mid\lambda''_{j}\mid^{\frac{\beta_{ij}}{3}}},\,\,\,\, b_{2}^{(i,j)}=\frac{\mid\lambda_{ij}\mid c_{3} c_{1}^{\frac{1}{2}}}{\mid\lambda''_{j}\mid^{\frac{1}{2}}} \,\,\,\,\hbox{and}\,\, \,\, b_{3}^{(i,j)}=\frac{\mid\lambda'_{ij}\mid c_{1}^{\frac{2}{3}}}{\mid\lambda''_{j}\mid^{\frac{2}{3}}},\,\, \hbox{for \,\,all}\,\, (i,j)\in \Omega_{n}.$$

Now, we state a straightforward, but useful result.
\begin{pro} \label{l1}
The Gribov  operator  matrix $\mathcal{M}^{\beta}_{\lambda'',\lambda',\lambda,\mu}$ is  diagonally generalized subordinate.\hf
\end{pro}
\begin{pf}
According to Proposition \ref{ff},  $H^{\beta_{ij}}_{\lambda'_{ij},\lambda_{ij},\mu_{ij}}$ is $\{\frac{\beta_{ij}}{3},\frac{1}{2},\frac{2}{3}\}$-generalized subordinate to $G_{\lambda''_{j}}$ for  all $(i,j)\in \Omega_{n}$. $\mathcal{M}^{\beta}_{\lambda'',\lambda',\lambda,\mu}$ is then diagonally generalized subordinate.
\end{pf}
\vskip 0.4 cm
As a consequence,  we obtain the following  result.
\begin{pro}\label{p4}
$\mathcal{R}^{\beta}_{\lambda',\lambda,\mu}$ is $\{\frac{\beta_{ij}}{3}\}_{(i,j)\in \Omega_{n}}\cup\{\frac{1}{2},\frac{2}{3}\}$-generalized subordinate to $\mathfrak{D}_{\lambda''}$ with bound $\displaystyle\{b_{1}^{(i,j)}\}_{(i,j)\in \Omega_{n}}\bigcup\{\displaystyle\sum_{(i,j)\in \Omega_{n}}b_{2}^{(i,j)},\displaystyle\sum_{(i,j)\in \Omega_{n}}b_{3}^{(i,j)}\}$\hf
\end{pro}
\begin{pf}
Follows immediately from   Lemma \ref{l0}, Proposition  \ref{l1} and the decomposition \eqref{GMAT}.
\end{pf}
\vskip 0.3 cm

\begin{rem}
Note that  0 is an eigenvalue of the operator matrix $\mathfrak{D}_{\lambda''}$  and that  $\mathcal{R}^{\beta}_{\lambda',\lambda,\mu}$ is generalized subordinate  but not $p$-subordinate to $\mathfrak{D}_{\lambda''}$. Therefore, the results of \cite{M1} and \cite{w}   are  not applicable  to the Gribov operator matrix $\mathcal{M}^{\beta}_{\lambda'',\lambda',\lambda,\mu}$.
\end{rem}

\vskip 0.3 cm

Below, we provide sufficient  conditions for the closedness of the Gribov  operator  matrix $\mathcal{M}^{\beta}_{\lambda'',\lambda',\lambda,\mu}$.
\begin{pro}\label{p5}
Suppose that $\displaystyle\sum_{(i,j)\in \Omega_{n}}\Big(\frac{1}{3}b_{1}^{(i,j)}\beta_{ij}+\frac{1}{2}b_{2}^{(i,j)}+\frac{2}{3}b_{3}^{(i,j)}\Big)<1.$ Then
the Gribov  operator  matrix $\mathcal{M}^{\beta}_{\lambda'',\lambda',\lambda,\mu}$ is closed.\hf
\end{pro}
\begin{pf}
 $\mathfrak{D}_{\lambda''}$ is  closed and $\mathcal{M}^{\beta}_{\lambda'',\lambda',\lambda,\mu}$ is diagonally generalized subordinate. So, it suffices to apply    \cite[IV Theorem 1.1]{TT}.
\end{pf}

\vskip 0.3 cm

The self-adjointness or non-self-adjointness of $\mathcal{M}^{\beta}_{\lambda'',\lambda',\lambda,\mu}$  depends on the non-symmetrical  operators $H_{\lambda_{ij}}:=\lambda_{ij}H_{1}$.  This has  of course  an effect on the location of the spectrum as well as on the existence of bases.
\begin{thm}\label{t7} The Gribov  operator  matrix $\mathcal{M}^{\beta}_{\lambda'',\lambda',\lambda,\mu}$ has compact resolvent. Moreover, one of the two following assertions holds: \begin{enumerate}
  \item If  for all $(i,j)\in \Omega_{n}$, $ \lambda_{ij}=0 $ and $\displaystyle\sum_{(i,j)\in \Omega_{n}}\Big(\frac{1}{3}b_{1}^{(i,j)}\beta_{ij}+\frac{2}{3}b_{3}^{(i,j)}\Big)<1,$    then $\mathcal{M}^{\beta}_{\lambda'',\lambda',0,\mu}$ is self adjoint and hence $\displaystyle\sigma(\mathcal{M}^{\beta}_{\lambda'',\lambda',\lambda,\mu})\subset\mathbb{R}$.
  \item   If  not for all $(i,j)\in \Omega_{n}$, $ \lambda_{ij}=0 $, then for every
$\alpha>\displaystyle\sum_{(i,j)\in \Omega_{n}}\Big(b_{1}^{(i,j)}+b_{2}^{(i,j)}+b_{3}^{(i,j)}\Big),$
there exists $r_0> 0$ satisfying\vskip 0.3 cm

\hfill$\displaystyle{\sigma(\mathcal{M}^{\beta}_{\lambda'',\lambda',\lambda,\mu})\subset B_{r_0}(0)\cup\{  (x+\textbf{i}y);\  x\geq0\mbox{ and } |y|\leq \alpha\max_{(i,j)\in \Omega_{n}} \{x^{\frac{2}{3}},x^{\beta_{ij}}\}\}.}$\hfill\hf
\end{enumerate}

\end{thm}
\begin{pf}
Firstly, since $\mathfrak{D}_{\lambda''}$ has   compact resolvent,  then $\mathcal{M}^{\beta}_{\lambda'',\lambda',\lambda,\mu}$ has so.
  For the first case, $\mathfrak{D}_{\lambda''}$ is self-adjoint and  $\mathcal{R}^{\beta}_{\lambda'',\lambda',0,\mu}$ is symmetric. Then the result follows from Kato-Rellich  Theorem \cite{Rellich}.
 For the second case, $\mathfrak{D}_{\lambda''}$ is self-adjoint with compact resolvent and  $\mathcal{R}^{\beta}_{\lambda'',\lambda',\lambda,\mu}$ is non-symmetrical. So, we should apply   Proposition \ref{p4} and  Theorem \ref{t4}.
\end{pf}

\vskip 0.4 cm

\begin{thm} \label{t8} One of the two following assertions holds:
\begin{enumerate}
  \item   If  for all $(i,j)\in \Omega_{n}$, $ \lambda_{ij}=0 $ and $\displaystyle\sum_{(i,j)\in \Omega_{n}}\Big(\frac{1}{3}b_{1}^{(i,j)}\beta_{ij}+\frac{2}{3}b_{3}^{(i,j)}\Big)<1$, then   $\mathcal{M}^{\beta}_{\lambda'',\lambda',0,\mu}$ admits  an orthonormal basis of eigenvectors.
  \item If  not for all $(i,j)\in \Omega_{n}$, $ \lambda_{ij}=0 $, then $\mathcal{M}^{\beta}_{\lambda'',\lambda',\lambda,\mu}$    admits a Riesz basis with parentheses.\hf\end{enumerate}
\end{thm}
\begin{pf}
  For the first case,  since $\mathcal{M}^{\beta}_{\lambda'',\lambda',0,\mu}$ is self-adjoint with compact resolvent, then we apply the spectral theorem.
 For the second case, let $1\leq j\leq n  $ and $\xi_k:=\lambda''_{j}(k-2)(k-1)k$  the eigenvalue number $k$ of $G_{\lambda''_{j}}$. Put
$r_k=\frac{\xi_k+\xi_{k+1}}{2}.$ We have

$\hfill\displaystyle\frac{N(r_k,G_{\lambda''_{j}})}{r_{k}^{1-\frac{2}{3}}}=\frac{2^{\frac{1}{3}}}{\lambda_{j}''^{\frac{1}{3}}}\frac{k}{({k(k-1)(2k-1)})^{1-\frac{2}{3}}}
\underset{r_k \longrightarrow
\infty}{\sim} \frac{1}{\lambda_{j}''^{\frac{1}{3}}}<+\infty.\hfill$

 Hence,
To finish,  it suffices to apply   Proposition \ref{p4} and  Theorem \ref{t5}.
\end{pf}
\vskip 0.4 cm
\begin{rem}
  By virtue of Theorems \ref{t7} and \ref{t8}, the operator matrices  $\mathfrak{D}_{\lambda''}$, $\mathfrak{S}_{\lambda'},$ $\mathfrak{H}_{\lambda}$ and $\mathfrak{H}_{\mu}^{\beta}$, have  the same spectral properties of their corresponding  entries.
\end{rem}

 \vskip 0.2 cm
We close this paper by a simple example of Gribov operator matrix by taking particular  values of the real  parameters    of Pomeron.
\begin{ex}
We consider the Gribov matrix:
$$\mathcal{M}^{3p}_{\lambda'',\lambda',0,\mu}= \lambda''Diag(G,\ldots,G)+\big((1-\delta_{ij})p_{ij}  H_0^{3p_{ij}}\big)_{1\leq i,j\leq n} ,$$
where  $\displaystyle\{p_{ij}\}_{(i,j)\in \Omega_{n}}$ is  a double indexed and decreased sequence defined by: $$\left\{\begin{array}{l} p_{12}=p_{21}=\frac{1}{3},\\\\
p_{ij}=\frac{1}{a^{i+j}},\ \mbox{for } (i, j)\in \Omega_{n}\setminus\{(1,2),(2,1)\}\,\,\, (a\geq\frac{7}{5}).\end{array}\right.$$
\end{ex}
\begin{pro}\label{p6}
Suppose that $\gamma:=\frac{c_{1}}{\mid\lambda''\mid}<  1$, then  $\sigma\big(\mathcal{M}^{3p}_{\lambda'',\lambda',0,\mu}\big)\subset \mathbb{R}$ and $\mathcal{M}^{3p}_{\lambda'',\lambda',0,\mu}$ admits an orthonormal basis of eigenvectors.\hf
\end{pro}
\begin{pf}
According to  Theorems \ref{t7} and \ref{t8},    it suffices to  prove that  $\small{\displaystyle\sum_{(i,j)\in \Omega_{n}}\gamma^{p_{ij}}p_{ij}^{2}<1.}$\\ \noindent\normalsize  We have
$$\sum_{(i,j)\in \Omega_{n}}\gamma^{p_{ij}}p_{ij}^{2}\leq2\sum_{1\leq i < j\leq n}p_{ij}^{2}=2p_{12}^{2}+2\sum_{2\leq i<j\leq n}p_{ij}^{2}=\frac{2}{9}+2\sum_{2\leq i < j\leq n}p_{ij}^{2}.
$$
\noindent Put  $S:=\displaystyle\sum_{2\leq i < j\leq n}p_{ij}^{2}$. Then
$S=\displaystyle\sum_{i=2}^{n-1}\frac{S_i}{a^{2i}}$, where $S_i=\displaystyle\sum_{j=i+1}^{n}\frac{1}{a^{2j}}$. Thus,
$$S_i=\frac{1}{a^{2(i+1})}\frac{1-(\frac{1}{a^{2}})^{n-i}}{1-\frac{1}{a^{2}}}<\frac{1}{a^{2i}(a^{2}-1)}.$$
Hence,
$$S< \frac{1}{a^{2}-1} \sum_{i=2}^{n-1}\frac{1}{a^{4i}}< \frac{1}{a^{4}(a^{4}-1)(a^2-1)}
.$$\noindent
Since $a\geq\frac{7}{5}$, then $S<\frac{7}{18}$ and therefore, $\small{\displaystyle\sum_{(i,j)\in \Omega_{n}}\gamma^{p_{ij}}p_{ij}^{2}<1.}$
\end{pf}

\end{document}